# Plug-n-Play Alternating Projection Algorithm for Large-scale Security Constraint Optimal Power Flow

Tuyen Vu, *Member, IEEE*

*Abstract*—In this paper, we present an optimization algorithm based on an alternating projection method to solve the large-scale security constraint optimal power flow (SCOPF) problem in power systems. The SCOPF is first partitioned into sub-problems, which share common power components. The proposed algorithm is fundamentally a distributed computing algorithm, which optimizes power flow in each electrical bus. The output of the local optimization problem will be exchanged among neighboring buses for the overall solution for the whole power system. We performed the algorithm against a modified 14-bus power system and demonstrated that the method is competitive with the commercial products regarding the CPU computation time.

*Index Terms*—SCOPF, ADMM, and projection algorithm.

## I. Introduction

### A. Overview

SCOPF has been a challenging problem for the power systems community as it is inherently nonlinear and non-convex optimization problem. In addition, the contingency analysis, such as N-1 contingency for the SCOPF normally requires a large number of iterations. Therefore, the computation time has been a challenge for scheduling power operation of generation systems in real-time. The objective of the work is to generalize the Projection Algorithm (PR) within the Alternating Direction Method of Multipliers (ADMM) framework to solve the SCOPF of large-scale power networks. In this paper, we propose to optimize the branch power flow variables instead of the regular generation variables. The SCOPF problem is partitioned into individual small problems, which represent the operational cost for each Bus. The technique optimizes branch power flows using PR in one Bus in exchanging their solution with neighboring Buses for the consensus power flow solution. Therefore, the technique is scalable to large-scale power networks. In addition, branch power flows are considered as the main variables for the problem, these variables can be effortlessly plugged on or switched off from the SCOPF solution. This plug-n-play approach allows the technique to be migrated to any power system models with minimal effort.

### B. State of the Art

In the current practices for the real-world systems, DCOPF has been the standard. The Independent System Operators (ISO) solve the DCOPF for the day-ahead market or real-time market [1]. However, DCOPF, as a linearized form of ACOPF does not provide the exact optimal solution. On the other hand, ACOPF is nonlinear and non-convex and very hard to solve. Recent advancements have been devoted for the ACOPF. The survey of OPF techniques can be found in [2]-[6]. The common techniques used are linear programming, quadratic programming, gradient-based methods, and Newton's method. Semidefinite programming (SDP) has been a promising method for this problem recently [7]-[9]. Fundamentally, SDP is based on the convex optimization and can be solved using interior point methods. For the large-scale problems, ADMM has been a promising method. The core of ADMM is in convex optimization and finding the Lagrangian multipliers [10]-[12]. However, the application of ADMM to the AC SCOPF is still limited.

In this paper, we will perform a new AC SCOPF reformulation based on the electrical line power flow components, which will be solved using an efficient Projection Algorithm based on an interior point method. The proposed technology will treat all Buses equally. The branch power flow is considered as the primary variable of the optimization; therefore, the method allows for seamless plug-n-play contingency analysis. The method also provides for straightforward incorporation of Distributed Energy Resources (DER) into the optimization problem. Hence, the method will apply to both transmission and distribution systems.

This paper is organized as follows: In Section II, the SCOPF problem will be described and partitioned into sub-problems, which can be solved separately in an effective manner with interior point method with projection algorithm to handle constraints. The ADMM is described next to coordinate the local solution for the global solution. In Section III, a 14-bus power system is taken for the case study with a demonstration that the method is competitive compared with the existing commercial products. Section IV will summarize our achievements.

## II. Proposed Method

Summary of the technique is shown in Fig. 1. In the figure, a three-bus system is taken as an example of the proposed approach. Each Bus is responsible for their own optimal cost. The solution of the cost function is the power flow among buses, which are exchanged among neighboring Buses. In the approach, the contingency tests for the SCOPF are done in parallel to identify constraints violation.

The cost function in each bus $i$ ($i = 1, 2, 3, ..., N$) without constraints is

$$c_i(p_i) = a_i p_i^2 + b_i p_i + c_i \tag{1}$$

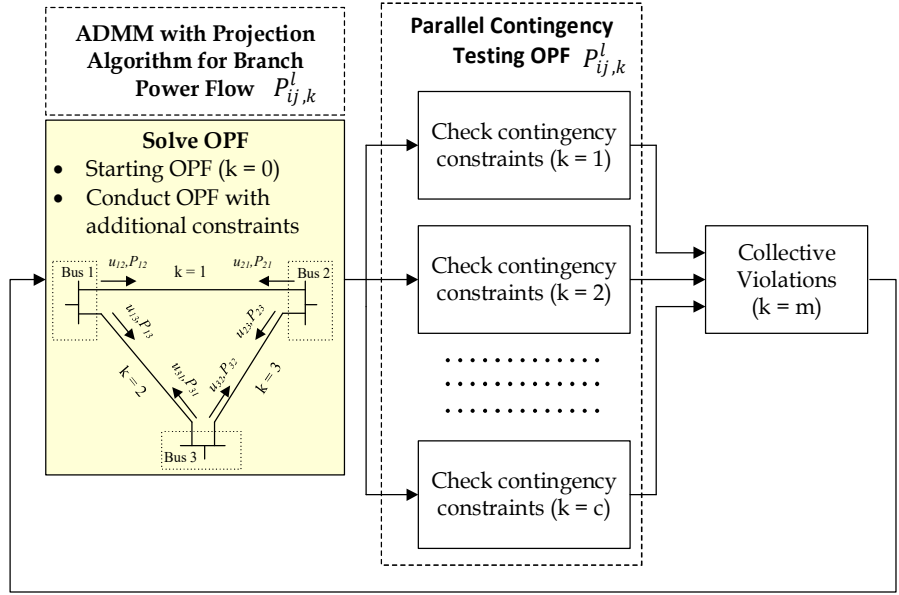

Fig. 1. Proposed SCOPF Solution Based on the Projection Algorithm within ADMM technique. A 3-Bus power system is an example for illustration.

Therefore, the formulation of the SCOPF problem can be expressed as

$$\min_{p_i} \left( c = \sum_{i \in G} c_i(p_i) \right)$$

$$\text{s.t.} \quad g_{o,k}(x_k, p_i) = 0, \quad k \in \{0\} \cup \mathcal{K}$$

$$h_{o,k}(x_k, p_i) \leq h_{o,k}^{max} \quad k \in \{0\} \cup \mathcal{K}, \tag{2}$$

where $i$ is Bus index, $\mathcal{K}$ is a set of contingency, $x_k$ represents Bus voltage magnitudes and angles under pre-contingency ($k = 0$) and contingency ($k \in \mathcal{K}$) conditions, $p_i$ is the power set under pre-contingency conditions, $c_i(p_i)$ is the cost of generation at Bus $i$, $g_{ok}(x_k, p_i) = 0$ represents the power flow, and $h_{ok}(x_k, p_i) \leq h_{ok}^{max}$ represents the branch limitations.

As contingencies are the interests for the solution, it is advantageous to include the branch variables into the optimization problem instead of generation variable. Fig. 1 shows the proposed technology for solving the problem. Rewrite (2) as a problem can be solved using ADMM.

$$\min_{p_{ji,k}, z_{ji,k}} \sum_{j \in N_i} \left( c_i(p_{ij,k}) + o_i(z_{ij,k}) \right.$$

$$\left. + \frac{\rho}{2} \| p_{ij,k} - z_{ij,k} + u_{ij,k} \|_2^2 \right) \tag{3a}$$

$$\text{s.t.} \quad p_{ij,k} \in \mathcal{C}_{i,k}$$

$$p_{ij,k} - z_{ji,k} = 0, \tag{3b}$$

where $\mathcal{C}_{i,k}$ represents $h_{i,k}(x_{i,k}, p_{ij,k}) \leq h_{ik}^{max}$, $p_{ij,k} - z_{ij,k} = 0$ represents the power flow constraint, and $o_i$ is the indicator function of $\mathcal{C}_i$, $\rho$ is a positive penalty number, $u_{ij,k} = (1/\rho)\lambda_{ij,k}$ is the scaled dual variable, and $N_i$ is the number of neighboring Buses. The solution is as

$$p_{ij,k,h+1} = \underset{p_{ji,k}}{\operatorname{argmin}} \left( c_i(p_{ij,k}) \right.$$

$$\left. + \frac{\rho}{2} \| p_{ij,k} - z_{ij,k,h} + u_{ij,k,h} \|_2^2 \right) \tag{4a}$$

$$z_{ij,k,h+1} = \underset{z_{ij,k} \in \mathcal{C}_{i,k}}{\operatorname{argmin}} \left( \sum_{j \in N_i} \left( o_i(z_{ij,k}) \right. \right.$$

$$\left. \left. + \frac{\rho}{2} \| p_{ij,k,h+1} - z_{ij,k} + u_{ij,k,h} \|_2^2 \right) \right) \tag{4b}$$

$$u_{ij,k,h+1} = u_{ij,k,h} + (p_{ij,k,h+1} - z_{ij,k,h+1}), \tag{4c}$$

where $h$ is the iteration index. It is important to note that the analytic solution of $z_{ij,k,h+1}$ using gradient method will result in an intermediate variable $z_{ij,k,h+1,0}$. To satisfy the contingency constraints set $\mathcal{C}_{i,k}$, $z_{ij,k,h+1,0}$ will be projected on to the set for a solution $z_{ij,k,h+1,1}$. In this paper, we perform a convex-relaxation for the $\mathcal{C}_{i,k}$; therefore, the solution $z_{ij,k,h+1,1}$ is unique.

$$\min \frac{1}{2} (z_{ij,k,h+1,0} - z_{ij,k,h+1,1}) \tag{5a}$$

$$\text{s.t.} \, z_{ij,k,h+1,1} \in \mathcal{C}_{i,k} \tag{5b}$$

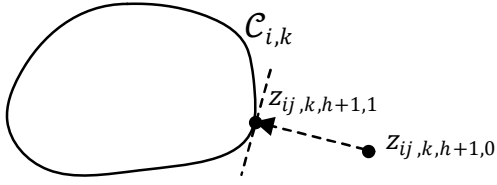

Fig. 2. Euclidean projection of $z_{ij,k,h+1,0}$ on $\mathcal{C}_{i,k}$.

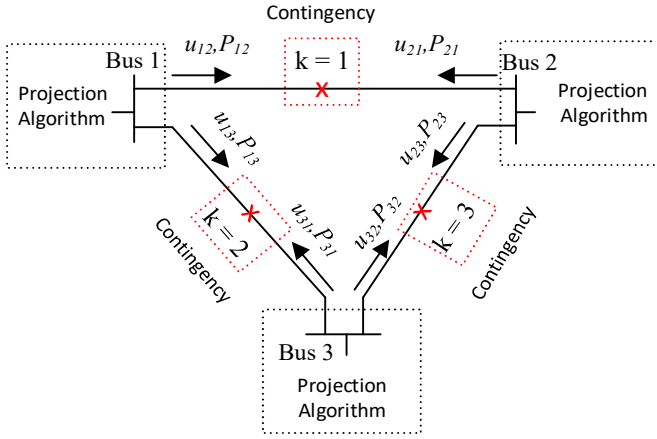

Fig. 3. An example of a 3-Bus power system for contingency.

Table 1. Step-by-step solution.

| Step | Optimization action at a Bus $i$ |
|---|---|
| 1. | Initialize parameters $z_{ij,k,h}$, and $u_{ij,k,h}$ (k=0) |
| 2. | Solve (4a) for $p_{ij,h+1}$, and $u_{ij,k,h}$ |
| 3. | Exchange $p_{ij,k,h+1}$, $u_{ij,k,h}$ with neighboring Bus $j$ |
| 4. | Solve (4b) for $z_{ij,k,h+1}$ |
| 5. | Update scale dual variable $u_{ij,k,h+1}$ |
| 6. | If the stopping criteria is satisfied go to step 7. Else, go to step 2 |
| 7. | Go back to perform steps 1-6 and run the algorithm for (k $\in \mathcal{K}$ ) and collect constraints if redispatch is required |
| 8. | Use the collected constraints to perform steps 1-7 |

Fig. 3 shows an example of a 3-Bus system. The solution variables $p_{ij,k}$ are generated in each Bus using (4a). The global variable $z_{ij,k}$ for consensus optimization is then calculated using (4a). Then the exchanged variable $z_{ij,k}$ is projected to the constraint set $\mathcal{C}_{i,k}$ for the feasible solution.

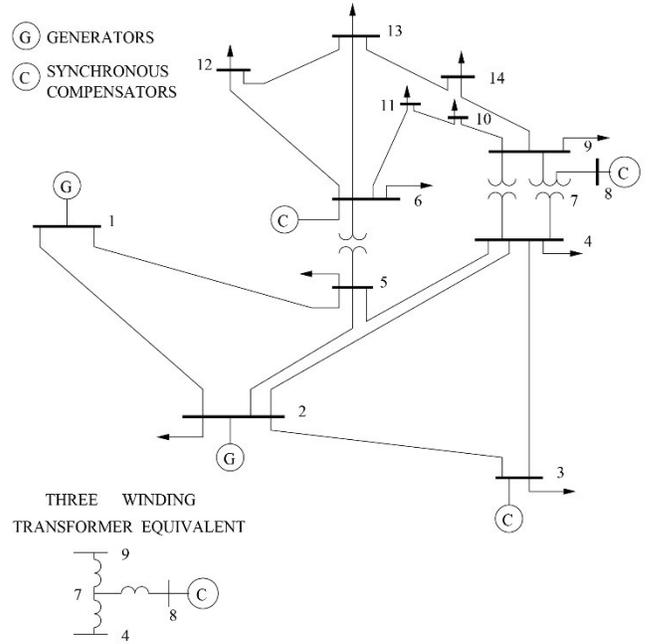

Fig. 4. IEEE 14-Bus.

Table 2. Correspondent of the cost and scenario for 14-Bus.

| Scenario | 1 | 20 | 40 | 60 | 80 | 100 |
|---|---|---|---|---|---|---|
| Cost ($) | 14,646 | 33,988 | 31,215 | 21,530 | 21,530 | 17,420 |

## III. CASE STUDY

To illustrate the effectiveness of the proposed approach, a modified 14-Bus system, posted for Beta phase challenge was considered for the SCOPF [16]. The DC OPF was performed. Data from scenarios 1, 20, 40, 60, 80, and 100 for the 14-Bus was utilized. A Dell laptop computer XPS 15 (Core i9-8950HK-2.9Ghz, 32GB RAM, 64-bit Windows) was used to verify the developed algorithm. In these scenarios, we performed two contingencies (one is the branch 4-5, and the other is the transformer 5-6). Cost of operation for these scenarios is shown in Table 2.

Regarding the capability, the minimum CPU time for the base case OPF was 0.03s and for the SCOPF (2 line contingencies) was 0.12s. Proof of the convergence for the optimal solution is shown in Fig. 5. In the figure, Residuals, $\|r_{i,h}\|_2^2$ and $\|s_{i,h}\|_2^2$ represent the power flow consensus among Bus and the convergence of the solution, respectively.

We performed the SCOPF on a Laptop computer for the 14-Bus system. To make a comparison between the technology and the existing commercial software, the CPU time for the base OPF is normalized as 0.03/14 = 2.1 (ms/Bus). The common commercial solvers are Knitro, Minos, Ipopt, and

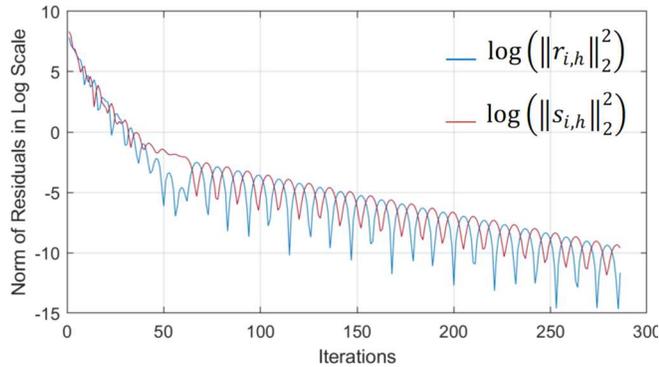

Fig. 5. Residual profiles in log scale of scenario 1 (Lower value is better).

Conopt. [15] performed OPF on these solvers on a Sever (Xeon E7458-2.4Ghz, 64GB RAM). The normalization was performed for the uniform starting points cases of these solvers. The normalized time for Knitro, Minos, Ipopt are 2.5, 1.9, and 1.5 (ms/Bus), respectively. Although produced by a Laptop computer, our results are comparable with these solvers. Therefore, we expect that the algorithm with Sever computer will be more competitive.

For the contingency cases, the CPU timing for 2 contingency cases was 0.12s. We anticipate that using a high-end Sever computer, the solution can be achieved within 5 minutes for a power network size of 100,000 Buses. The uniqueness of the proposed method is that it allows for plug-n-play condition. The method is generalized for every single Bus in the system. Therefore, adding or removing one-line for SCOPF can be done at a minimum effort.

As illustrated the proposed algorithm is effective and competitive regarding CPU time compared with the existing commercial tools.

## IV. Conclusion

As the proposed technology is scalable, the work is valuable not only for the transmission systems but also to the distribution systems. The distribution systems increasingly incorporate distributed energy resources such as energy storages and distributed resources are interconnected. In this case, millions of devices will need to be optimized. Therefore, this work will also be valuable for the distribution systems.